\documentclass[11pt]{article}
\usepackage{url}
\usepackage{calc}
\usepackage{a4wide}
\usepackage{listings}
\usepackage[ngerman,english]{babel}
\usepackage{etex}
\usepackage{ellipsis,ragged2e}
\usepackage{microtype}
\usepackage{amsmath}
\usepackage{amssymb}
\usepackage{amsthm}
\usepackage{color}
\usepackage{times}
\usepackage{graphicx}
\usepackage{caption}
\usepackage{hyperref}
\usepackage{cite}
\usepackage{mathrsfs,graphicx,color,float, indentfirst,textcomp}
\usepackage{setspace}
\usepackage{latexsym,amssymb,lscape,amsmath,amsthm,amsfonts,amssymb,cite,enumerate}
\usepackage{lscape}
\usepackage{xfrac}
\usepackage{indentfirst}
\usepackage[caption=false]{subfig} 

  \usepackage{tikz}
  \usetikzlibrary{arrows}
  \usepackage{pgfplots}
  \pgfplotsset{compat=newest}
  \pgfplotsset{plot coordinates/math parser=false}
  \newlength\figureheight
  \newlength\figurewidth
  
\newtheorem{thm}{Theorem}[section]

\newtheorem{lem}[thm]{Lemma}

\newtheorem{defn}[thm]{Definition}

\newtheorem{example}[thm]{Example}

\newcommand{\R}{\mathbb{R}}
\renewcommand{\u}{{\bf U}}
\newcommand{\uad}{{\mathcal{U} }}
\renewcommand{\v}{{\bf V}}
\newcommand{\x}{{\bf x}}
\newcommand{\eps}{{\varepsilon}}
\newcommand{\srw}{{\delta}}

\newcommand{\us}{ ({\bf U^\pm_x})  } 
\newcommand{\ar}[2]{\left(\begin{array}{#1}#2\end{array}\right)}
\newcommand{\Mach}{{Ma}} 
\newcommand{\gammaMinusOne}{\gamma-1} 
\newcommand{\fig}[1]{Figure~#1}

\def\figscale{0.8}

\begin{document}

\title{Algorithmic differentiation of hyperbolic flow problems}
\date{\today}	
\author{Michael Herty\thanks{RWTH Aachen University, Institute of Geometry and Applied Mathematics, 
		Templergraben 55, 52056 Aachen, Germany, herty@igpm.rwth-aachen.de}
	\and 
Jonathan H\"user \thanks{RWTH Aachen University, Informatik 12: Software and Tools for Computational Engineering (STCE), 
	52056 Aachen, Germany, hueser@stce.rwth-aachen.de}
\and 
Uwe Naumann \thanks{RWTH Aachen University, Informatik 12: Software and Tools for Computational Engineering (STCE), 
	52056 Aachen, Germany, naumann@stce.rwth-aachen.de}
\and 
Thomas Schilden \thanks{RWTH Aachen University, Chair of Fluid Mechanics and Institute of Aerodynamics, 
	52062 Aachen, Germany, t.schilden@aia.rwth-aachen.de}
\and 
Wolfgang Schr\"oder \thanks{RWTH Aachen University, Chair of Fluid Mechanics and Institute of Aerodynamics, 
	52062 Aachen, Germany and RWTH Aachen University, JARA Center for Simulation and Data Science, 52074 Aachen, Germany, office@aia.rwth-aachen.de}
\and 
}  

\maketitle

\begin{abstract}
We are interested in the development of an algorithmic differentiation  framework for computing approximations to  tangent vectors to scalar and systems of hyperbolic partial differential equations.   The main difficulty of such a numerical method is the presence of shock waves that are resolved by proposing a numerical discretization of the calculus introduced in Bressan and Marson~[Rend. Sem. Mat. Univ. Padova, 94:79--94, 1995].
Numerical results are presented for the one-dimensional Burgers equation and the Euler equations.  Using the essential routines of a state-of-the-art code for computational fluid dynamics (CFD) as a starting point, three modifications are required to apply the introduced calculus.
First, the CFD code is modified to solve an additional equation for the shock location. Second, we customize the computation of the corresponding  tangent to the shock location. Finally, the modified method is enhanced by algorithmic differentiation.  
Applying the introduced calculus to problems of the Burgers equation and the Euler equations, it is found that correct sensitivities can be computed, whereas the application of black-box algorithmic differentiation fails.
\end{abstract}
	

\noindent {\bf AMS}
  35L65, 49K20, 49K40

\section{Introduction}
We are interested in an algorithmic differentiation framework for the computation of sensitivities to multi-dimensional systems of hyperbolic partial differential equations. Such a framework is relevant, for instance, in supersonic flows which are often characterized by the occurrence of shock waves. 
Prominent examples are external flows over trans- to supersonic aircraft and internal supersonic flows through nozzles or diffusers of, e.g., ramjets.
Given a design parameter for an arbitrary objective function, the sensitivity of the objective function with respect to the design parameter has to account for the discontinuities introduced by the shocks.
 
Towards developing such a method  we are concerned in this work with a suitable  algorithmic framework for scalar but possible multi-dimensional hyperbolic problems. The prototype of this problem is defined by equation \eqref{claw} where we denote by $u(t,x) \in \R$  the unique entropy solution.  The flux $f \in C^4(\R^d; \R^d)$ is assumed to be nonlinear. 
\begin{equation}
\begin{aligned} 
\partial_t u + \nabla_x \cdot f(u)&=0,&&x \in \mathbb{R}^d,\quad t>0,\\
u(0,x)&=u_0(x),&&x \in \mathbb{R}^d.
\end{aligned}
\label{claw}
\end{equation}

In the spatially one-dimensional case $(d=1)$ there has been tremendous progress in both analytical and numerical studies of problems of sensitivities of $u$ with respect to 
initial data $u_0,$ see e.g.,  \cite{BandaHerty2012aa,BardosPironneau2002aa,Bianchini2000aa,BressanGuerra1995aa,BressanLewicka1999aa,BressanMarson1995ab,BressanMarson1995aa,CastroPalaciosZuazua2008aa,GilesUlbrich2010aa,
	JamesSepulveda1999aa,Ulbrich2003aa,PierceGiles2004aa,Giles1996aa,GilesSueli2002aa,GugatHertyKlar2006aa}.
Even in the one-dimensional, scalar case, it has been shown that  the evolution operator $\mathcal{S}_t: u_0(\cdot)\to u(\cdot,t)= \mathcal{S}_t u_0(\cdot)$ generated by  the conservation law is generically non--differentiable in $L^1$ \cite[Example 1]{BressanMarson1995ab}.  A theoretical calculus for
the first-order sensitivities of $\mathcal{S}_t u_0$ with respect to $u_0$ has been established in \cite[Theorems 2.2 and 2.3]{BressanMarson1995ab} for
general spatially one-dimensional systems of conservation laws. Here, the initial data $u_0$ is assumed to be  piecewise Lipschitz continuous 
and contains finitely many discontinuities. Therein, the concept of  tangent vectors has been introduced to characterize the evolution of variations with respect to
$u_0,$ see \cite[equations (2.16)--(2.18)]{BressanMarson1995ab}. It has been further extended in \cite{BressanCrastaPiccoli2000aa}
to establish  continuous dependence of $\mathcal{S}_t$ on the initial data $u_0.$ This result has been extended to BV (bounded variation) initial data in
\cite{BressanGuerra1995aa,Bianchini2000aa} and led to the introduction of a differential structure for $u_0\to \mathcal{S}_t u_0$, called
shift-differentiability, see e.g. \cite[Definition 5.1]{Bianchini2000aa} and 
an adjoint calculus  \cite[Proposition 4]{BressanShen2007aa}. In the scalar, one-dimensional 
case the assumptions on $u_0$ could be weakened as shown e.g. in \cite{Ulbrich2003aa,CastroPalaciosZuazua2008aa}.  Analytical results for optimal control problems in the 
case of a one-dimensional, scalar hyperbolic balance laws with a convex flux have also been developed using a different approach in \cite{Ulbrich2003aa}. The relation to the weak formulation has been discussed in \cite{BardosPironneau2002aa} for the Burgers equation. 
\par 
The theoretical sensitivity calculus provides equations for the evolution of the variation of the value of the solution 
coupled to the evolution of the  variation of the positions of possible shocks in $u$. This provides evolution equations for
the  tangent vector $(v,\xi)$ of $S_t u_0$ at $u_0$. We present here an algorithmic framework that allows for a {\em numerical computation } using
algorithmic differentiation. This requires in  particular,   to augment a possible numerical simulation code for equation \eqref{claw}  by an  evolution for the  possible shock positions. This also requires to change the notion of forward differentiability.  The augmentations will be described in the  algorithmic differentiation framework introduced below. 
\par 
Black--box algorithmic differentiation (AD) \cite{Griewank2008EDP,Naumann2012TAo} assumes (classical) differentiability 
of the mapping $u_0 \to S_t u_0$ which is the main reason for it not being applicable to the given
scenario prior to the proposed modification. In the following we introduce the notation used in the AD framework.  Consider a nonlinear, finite--dimensional 
map $\uad \to G(\uad):\R^m \to \R^n$ and denote by $\uad_m=G(\uad_0).$ A
 directional derivative of $\uad_m$ with respect to  $\uad_0$  will then be denoted by 
$$ \dot{\uad}_m = \frac{d}{d \uad} G(\uad_0) .$$
Clearly, there is also an adjoint formulation that is preferable if gradients of scalar control objectives
are required.  AD has been applied successfully to numerous real-world applications in computational science, engineering and finance; refer, e.g., to  \cite{Bischof2008AiA,Bucker2005ADA,Forth2012RAi} for
further reference.  Software tools for AD use either source code transformation, e.g, 
\cite{Hascoet2013TTA} or function and operator overloading if supported by the 
programming language, e.g, \cite{Griewank1996AAC}. 
The computational experiments reported in Section~\ref{sec:cr} are based on
dco/c++ \cite{AIB-2016-08}.
Collections of both applications of AD and of AD software tools can be found on the community's web 
portal {\tt www.autodiff.org} together with a comprehensive bibliography on 
the subject.
\par 
Direct numerical methods based on the discretization of the tangent equations 
have been discussed e.g. in \cite{BandaHerty2012aa,HertyPiccoli2016aa,Giles1996aa,CastroPalaciosZuazua2008aa}. In \cite{GilesUlbrich2010aa}, the adjoint equation
has been discretized using a Lax-Friedrichs-type scheme, obtained by including conditions along shocks and modifying the Lax-Friedrichs
numerical viscosity.  Convergence results for sensitivity and adjoint equations have been obtained in \cite{Ulbrich2003aa} for a general class of finite--volume schemes satisfying a one-sided Lipschitz condition (OSLC) and in \cite{BandaHerty2012aa,HertyPiccoli2016aa}  
for implicit-explicit finite-volume methods.
Other examples of finite volume methods and Lagrangian methods are given e.g. in  \cite{ChertockHertyKurganov2014aa,HERTYKURGANOVKUROCHKIN2014aa}.
Compared to previous methods, we aim to suitably augment a forward simulation of a standard finite--volume discretization of equation \eqref{claw}, such that AD  yields  tangent vectors consistent with the theoretical calculus proposed in \cite{BressanMarson1995aa}. 

The paper is organized as follows. First, the concept of the new calculus is outlined. Then, the numerical implementation and software tools are given. The results of the numerical simulations are presented in \S~\ref{sec:cr}, before the findings are concluded.

\section{Theoretical Calculus}\label{sec::calculus}
We briefly recall the theoretical calculus and introduce the basic notion of tangent vectors. The presentation of the latter follows closely \cite{BressanMarson1995ab} and \cite{BressanCrastaPiccoli2000aa}. Then, theoretical results in dimension $d=1$ applied to scalar hyperbolic equations are derived and a numerical scheme is proposed. To increase the readability, we support the theory by applying the calculus to a prominent example of the Burgers equation~\cite{BressanMarson1995ab}. 
\begin{example} Here and in the first numerical results in Section \ref{susec::cr::burgers} we consider Burgers equation with $f(u)=\frac12 u^2$ and 
	a function $u_0(x)$ having  a single discontinuity $N(u_0)=1$ at $x_1=1$ given by 
	\begin{align}\label{eq::ramp::ic}
	u_0(x)= x \; \chi_{[0,1]}(x).
	\end{align}
	In a single spatial dimension a weak solution $u$ to \eqref{claw} that is in BV is a composition of piecewise Lipschitz continuous parts separated by jump discontinuities. Therefore in the following we consider $u_0$ of this class $\mathcal{U}.$  The particular structure of $u_0$ (and the corresponding solution $u$) suggests to consider variations of the Lipschitz parts as well as variations of the jump discontinuities. This motivates the  notion of a tangent space $T_u$ defined below and norm given by equation \eqref{norm}. The elements of this tangent space are called (generalized)  tangent vectors $(v,\xi) \in T_u$. In the following, we  are interested in applying AD for computing  a numerical approximation to 
	$(v,\xi).$  
\end{example}
Assume $f \in C^4(\mathbb{R})$ and $$u_0 \in \mathcal{U}:=\{ u: \mathbb{R} \to \mathbb{R}: u \mbox{ measurable }, TV(u)\leq C,   u \mbox{ piecewise Lipschitz continuous} \} ,$$
where $TV$ denotes the total variation.
  For  $u_0 \in \mathcal{U}$ we indicate by $x_k=x_k(u_0), k=1,\dots, N(u_0)$ the points of discontinuity of the function $u_0.$ 
For a function $u_0 \in \mathcal{U}$ a generalized tangent vector consists of 
two components $(v,\xi)$ where $v\in L^1(\mathbb{R})$ describes the $L^1$ infinitesimal displacement of $u_0$. Further,  $\xi \in \mathbb{R}^{N(u_0)}$ describes the infinitesimal displacement 
of the $N(u_0)$ discontinuities.  A norm on the space of tangent vectors $T_u:=L^1(\mathbb{R}; \mathbb{R}^n) \times \mathbb{R}^{N(u_0)}$ is given by 
\begin{align}\label{norm}
\| (v,\xi) \| := \| v \|_{L^1} + \sum\limits_{i=1}^{N(u_0)} | \Delta_i u_0 | \; | \xi_i |
\end{align}
where $\Delta_i u_0=u_0( x_i+) - u_0(x_i-)$.   The norm depends on $u_0$ through the number of points of discontinuities. 
\par 
Let $\eps>0$ be a sufficiently small parameter and let $(v,\xi) \in T_u$. Then, variations at $u_0$ are  
described by  shifting the function values by $ \eps  v$ and the $i$th discontinuity
by $   \eps   \xi_i.$ The resulting function $u^\eps$ is defined by 
\begin{align}\label{update}
u^\eps = u_0 + \eps \cdot   v -
\sum\limits_{i: \xi_i>0}^{N(u_0)} \Delta_i u_0 \; \chi_{ [x_i, x_i + \eps \xi_i ]} + \sum\limits_{i: \xi_i<0}^{N(u_0)} \Delta_i u_0\; \chi_{ [x_i + \eps \xi_i, x_i ]}.
\end{align} 
\begin{example}\label{ex::burgers::ramp} We may consider $\xi_1=0$ and $v(x)=x \chi_{[0,1]}(x).$ Then, the resulting function $u^\eps$ is given by 
	\begin{align}\label{u delta}
	u^\eps(x)= (1+\eps ) \; x \; \chi_{[0,1]}(x).
	\end{align}
\end{example}

   For $\eps$ sufficiently
  small $u^\eps$ has the same number of discontinuities as $u_0.$ Note that if $\xi \not = 0$ then the function $\eps \to u^\eps$
  is {\em not } differentiable in $L^1$ as stated in the introduction. In fact, the ratio  $\frac{ u^{\eps+ h} - u^\eps}h$ does not converge to any limit in $L^1$ for $h \to 0.$ 
  However, the previous limit remains meaningful  as a weak limit in a space of measures with a singular point mass located
  at $x_i$ and having magnitude $ | \Delta_i u | \xi_i.$ Therefore, in \cite{BressanMarson1995ab} a class of variations $\eps \to u_\eps$ is described up to 
  first order by  (generalized) tangent vectors $(v,\xi) \in T_u.$ 
   
  \begin{example}
  	Consider $u^\eps(x)$ as in equation \eqref{u delta}.	If we consider Burgers equation 
  	\begin{align}\label{eq::burgers}
  	\partial_t u + \partial_x \frac12 u^2 = 0, \; u(0,x) = u^\eps(x) 
  	\end{align}
  	we obtain for $t\geq0, x \in \R$ an explicit solution denoted by  $U^\eps(t,x)$ as 
  	\begin{align}\label{eq::burgers::solution_u}
  	U^\eps(t,x) = \frac{ (1+\eps) x}{1+ (1+\eps)t} \chi_{ [0, \sqrt{1+ (1+\eps)t} ]}(x).
  	\end{align} In fact, for $\eps \to U^\eps(0,x)\equiv u^\eps(x)$ is differentiable in $L^1$ and its derivative is precisely $v(x)$.   	
  	  However, for any positive $t>0$ the function $\eps \to U^\eps(t,\cdot)$ is {\em not} differentiable in $L^1(\R).$ The location of
  	  the discontinuity of $x\to U^\eps(t,x)$ and $x \to U^0(t,x)$ are different. The next paragraph discusses in which sense $U^\eps$ 
  	  can be expanded in terms of $\eps$ to allow for a characterization of the tangent $(v,\xi)$ of $U^\eps(t,\cdot)$ for $t>0.$
  	  The characterization is then given by equation \eqref{def tv 2}.
  \end{example}

  Let $u \in L^1(\mathbb{R}; \mathbb{R}^n)$ be a piecewise Lipschitz continuous function with $N=N(u)$  discontinuities. 
  Consider $\Sigma_u$, the family of all continuous paths $\gamma: [0,\eps_0] \to L^1_{loc}$ with $\gamma(0)=u$ with $\eps_0$ possibly
  depending on $\gamma.$ We recall  \cite[Definition 1,3]{BressanMarson1995ab}.
\begin{defn} 
	The space of generalized tangent vectors to a piecewise Lipschitz function  $u$ with jumps located at the points $x_1 < x_2 \dots < x_N$ is $T_u:=L^1(\mathbb{R}; \mathbb{R}^n) \times \mathbb{R}^{N}$. 
	A continuous path $\gamma \in \Sigma_u$  generates a tangent vector $(v,\xi) \in T_u$ if 
	$$ \lim\limits_{\eps \to 0} \frac{1}\eps  \| \gamma(\eps)  - \bar{\gamma}(\eps) \|_{L^1} =0$$
	for 
	\begin{align}\label{def tv}
	 \bar{\gamma}(\eps) := u + \eps v - \sum\limits_{i: \xi_i>0} \Delta_i u \; \chi_{ [x_i, x_i + \eps \xi_i ]} + \sum\limits_{i: \xi_i<0} \Delta_i u\; \chi_{ [x_i + \eps \xi_i, x_i ]}.
	 \end{align}
		\par 
	Let $u$ be a piecewise Lipschitz function with simple discontinuities \cite[Definition 2]{BressanMarson1995ab}. Then, a
	path $\gamma \in \Sigma_u$ is a regular variation for $u$ if additionally all function $\gamma(\eps)=u^\eps$ are piecewise Lipschitz with simple
	discontinuities   and the location of the jumps at  $x_i^\eps$ depend
	continuously on $\eps.$  
\end{defn}
A regular variation $\gamma$ for $u$ generates a tangent vector $(v,\xi)$  by 
\begin{align}\label{def tv 2}
\xi_i = \lim\limits_{\eps \to 0} \frac{ x_i^\eps - x_i }\eps, \; \lim\limits_{\eps \to 0} \int_{a}^b  
\left\| \frac{u^\eps(x_i^\eps + y) - u(x_i+y) }{ \eps } - v(x_i+y) - \xi_i u_x(x_i+y) \right\| dy = 0
\end{align}
whenever $[x_i+a,x_i+b]$ does not contain any other point of discontinuity of $u$ except $x_i.$ Further, the length of a regular path $\gamma$ can be computed 
by \eqref{norm}. We now consider the initial data $u_0$ and a regular variation generating the
tangent vector $(v,\xi) \in T_u.$

 \begin{example}\label{ex::burgers::ramp::solution}
	For $U^\eps(0,x)=u^\eps(x)$ the pair is $(v,\xi_1)$ where $v(x)=x \chi_{[0,1]}(x)$ and $\xi_1=0$ is a tangent vector by definition of $u^\eps.$ 
	Consider now $t>0.$ The position of the shock $x^\eps(t)$ in $U^\eps$  and the position $x(t)$ of the shock in the solution to Burgers' equation with initial datum  $u_0(x)$ are given by 
	\begin{align}\label{eq::burgers::shock_location}
	x^\eps(t) = \sqrt{1+ (1+\eps)t}  \mbox{ and } x(t) = \sqrt{1+ t}, 
	\end{align}
	respectively. 	Hence, the first term in equation \eqref{def tv 2} yields $\xi_1=\xi_1(t)$  
	\begin{align}\label{xi1}
	\xi_1(t) = \frac{t}{ 2 \sqrt{1+t} 
	}.
	\end{align}
	Furthermore, the second term in equation \eqref{def tv 2} yields $v(x)=v(t,x)$ as 
	\begin{align}\label{v1}
	v(t,x) =\frac{x}{(1 +  t)^2}  \chi_{[0,  \sqrt{1+t}]}(x). 
	\end{align}
The pair $(v,\xi_1)$ is the tangent vector to $u_0.$ It is computed using the explicit solution for the initial variation $(v(0,x)=x \chi_{[0,1]}(x),\xi_1(0)=0)$ introduced above. 
Lemma \ref{thm} shows that the tangent vector \eqref{xi1} and \eqref{v1} can also be obtained by propagating the initial   variation $\left(v(0,x),\xi_1(0) \right).$ This also yields a recipe for the AD tool: we might implement a suitable discretization, denoted by $G=(G_1,G_2)$ for the evolution of $t\to x(t)$ as well as for $u_0(\cdot) \to u(t,\cdot)$. Then, the directional derivatives of $x(\cdot)$ and $u(t,\cdot)$ with respect to $u_0$ lead to approximations of $(v,\xi_1).$ 
\end{example}

 Under  regularity assumptions  the regular variations $\gamma$  are locally preserved and  linearized equations exist for the  evolution of the tangent vector $t \to \left( v(t,\cdot), \xi(\cdot) \right)$. The following Lemma \ref{thm} is a consequence of \cite[Theorem 2.2]{BressanMarson1995aa}. 

\begin{lem}\label{thm}
   Consider equation \eqref{claw} for $d=1$ and $u(t,x)\in \R.$ 	Let $u(\cdot,\cdot)$ be a piecewise Lipschitz continuous solution to \eqref{claw}  and initial data  $u(0,x)= u_0(x)$ piecewise Lipschitz with $N=N(u_0)$
	simple discontinuities. Let $( \bar{v}, \bar{\xi} ) \in T_{u_0}$ be a tangent vector to $u_0$ generated by the regular variation $\gamma$ with $\gamma(\eps)= u_0^\eps$. Let $u^\eps(t,x)$ be the solution to \eqref{claw}
	and  initial data  $u^\eps(0,x) = u_0^\eps(x).$ Then, there exists a time $t_0>0$ such that for all $t \in [0,t_0]$  the path  $\bar{\gamma}$
	with $\bar{\gamma}(\eps) =   u^\eps(t, \cdot)$
	is a regular
	variation of $u(t,\cdot)$ generating the tangent vector $(v(t), \xi(t)) \in T_{u(t,\cdot)}$. Further,  $(v,\xi)$  is the unique (broad) solution to 
	\begin{align}
	\label{tv-v}
	v(0,\cdot) = \bar{v}(\cdot), \quad v_t + f(u) v_x + ( \frac{d}{du} f(u)v)u_x =  0, 
	\end{align} 
	outside of the discontinuities of $u$. For $i=1,\dots,N$ we have 
	\begin{align}\label{tv-xi}
	\frac{d}{dt}\xi_i(t) = \partial_{u^+} A(u^+,u^-) \left( v^+ + \xi_i(t) u^+_x\right) +   
	\partial_{u^-}  A(u^+,u^-) \left( v^- + \xi_i(t) u^-_x\right) 
	\end{align}
	along each line of discontinuity $x_i(t)$ where $u$ has a discontinuity. Here, $\Delta_i v = v^+-v^-$,  $v^\pm = v(x_i(t) \pm, t)$ and  	$A(u,v)=\int_0^1 \frac{d}{du} f( \theta u + (1-\theta)v) d\theta.$
\end{lem}
We refer to Definition \ref{broad solution} for the notion of broad solutions.  
Compared to  the general result \cite[Theorem 2.2]{BressanMarson1995aa} we note the following:  Due to fact that $u(t,x)\in \R$ we do not have discontinuities of different families and 
the consistency condition on the eigenvectors of $\frac{d}{du} f(u)$ is trivially satisfied. The full result
is given in the Appendix \ref{appendix} for convenience. 

\begin{example} Since $A(u,v)=\frac12 ( u +v )$  we have $\partial_{u}A= \partial_v A=\frac12.$  The shock position is $x(t)= \sqrt{1+t}$ and 
	the solution $u(t,x) = \frac{x}{2+t} \chi_{[0,  \sqrt{1+t}]}(x)$ and 
	therefore $u^+ = 0$ and $u^-=\frac{ \sqrt{1+t} }{ 2+ t}.$ Furthermore, $u^+_x=0$ and $u^-_x = \frac{1}{2+t}.$ 
	Hence, the corresponding equations for the example read
	\begin{align}
	v(0,x) = x \chi_{[0,1]}(x), \; v_t + \frac12 u^2 v_x + ( u v ) u_x = 0,  \label{v2} \\
	\xi_1(0)=0, \; \frac{d}{dt} \xi_1(t) =  \frac12 \left( v^+ +  v^- + \xi_1(t)  \frac{\sqrt{1+t} }{2+t} \right). \label{xi2}
	\end{align}
	One checks  that $v(t,x)$ given by \eqref{v1} fulfills equation \eqref{v2} pointwise except along $(t,x(t)).$ Since $v^+=0$ and $v^-=\frac{\sqrt{1+t}}{ (1+t)^2 }$  the right-hand side of 
	equation \eqref{xi2} is given by
	\begin{align}
	 \frac{1}{2 (1+t)^{\frac32}} +  \frac{t}{4(1+t)^{\frac32}} = \frac{d}{dt} \left( \frac{t}{2 \sqrt{1+t}}\right).
	\end{align}
\par 
The purpose of the AD framework applied in Section \ref{ad sec} is to avoid explicitly implementing equation \eqref{tv-v} and equation \eqref{tv-xi}. Equation \eqref{tv-v} is formally obtained by  linearizing equation  \eqref{claw}. Therefore, we expect that after suitable definition of the function $G_1$ black-box AD provides a suitable approximation, see e.g. equation \eqref{uwe_ad1}. However, in order
to obtain equation \eqref{tv-xi} we introduce  Lemma \ref{thm2} that shows a possible derivation of equation \eqref{tv-xi}. Those equations will lead to a further component $G_2$ of the numerical discretization $G.$
\end{example}

We observe that \eqref{tv-v} is the linearization of the forward dynamics \eqref{claw}. Therefore, we expect that if a  finite-volume scheme resolves the dynamics \eqref{claw} 
with sufficiently high accuracy, denoted by $G_1$, an AD tool  will produce a solution to equation \eqref{tv-v} with sufficiently high accuracy. However, as seen in Lemma \ref{thm} this only describes one component of the sensitivity $v$. Hence, we need to augment the AD by including equation \eqref{tv-xi} leading to a second component $G_2.$ Since we do not want to discretize equation \eqref{tv-xi} a posteriori we augment the forward simulation code by an additional computation of the shock location $x_i(t).$ In fact, the following Lemma holds true (see \cite{CastroPalaciosZuazua2008aa}).
  \begin{lem}\label{thm2}
  	Consider equation \eqref{claw} for $d=1$ and $u(t,x)\in \R.$ Assume that the function 
  	$u(t,\cdot) \in \mathcal{U}$ has a discontinuity at $x=x(t)$ and across the discontinuity 
  	the Rankine-Hugenoit condition is fulfilled:
  	\begin{align}\label{RH}
  	\frac{d}{dt} x(t) =\frac{ \Delta f(u) }{ \Delta u } 
  	\end{align}
  	where $\Delta u=u^+(x(t)+,t) - u^-(x(t)-,t).$  
  	Consider a regular variation $\gamma(\eps)=u^\eps$ of $u$ defined by equation \eqref{update}  with $N(u)=1$ and 	tangent vector $(v,\xi)$ and assume $u^\eps$ fulfills 
  	\begin{align}\label{RH-2}
  	\frac{d}{dt} x^\eps(t) =\frac{ \Delta f(u^\eps) }{ \Delta u^\eps } 
  	\end{align}
   Then, the first--order expansion in terms of $\eps$ of equation \eqref{RH-2}  is equivalent to equation \eqref{tv-xi}, i.e., 
 \begin{align}\label{tv-xi2}
 \frac{d}{dt}\xi(t) = \partial_{u^+} A(u^+,u^-) \left( v^+ + \xi(t) u^+_x\right) +   
 \partial_{u^-}  A(u^+,u^-) \left( v^- + \xi(t) u^-_x\right) 
 \end{align}
  	\end{lem}
  \noindent {\bf Proof.}
 Consider equation \eqref{tv-xi2}. According to the definition $A(\cdot,\cdot)$ we have for $u\not =v$ 
 \begin{align}
 \partial_{u}  A(u,v) =  
  \int_0^1 f''( \theta u + (1-\theta)v)  \theta d\theta = - \int_0^1 f'(\theta u + (1-\theta) v ) d\theta \frac1{u-v} +  \frac{f'(u)}{ u-v }, \\
  \partial_{v}  A(u,v) = \int_0^1 f''( \theta u + (1-\theta)v) (1-\theta) d\theta
  =  \int_0^1 f'(\theta u + (1-\theta) v ) d\theta \frac1{u-v} -  \frac{f'(v)}{ u-v }. 
 \end{align} 
  Hence, equation \eqref{tv-xi2} is equivalent to 
  \begin{align}
 \frac{d}{dt} \xi_i(t) =  
  \frac{  \Delta \left( f'(u) \left(   v +\xi_t u_x \right) \right)  }{ \Delta u }
  - \frac{ \Delta f(u) }{ (\Delta u)^2 } \Delta (v+\xi_t u_x) . 
  \end{align}
  Next, consider a regular variation $u^\eps$ of $u.$ For $\xi(t)\geq 0$ it is given by 
  \begin{align}\label{eq::calculus}
  u^\eps(t,x)= u(t,x) + \eps v(t,x) + \chi_{[x(t),x(t)+\eps \xi(t)]}(x) \Delta u
  \end{align}
  Due to the definition of $\xi$ we have 
  \begin{align}
  \lim\limits_{\eps \to 0} \frac{ x^\eps(t)-x(t) }\eps = \xi(t).
  \end{align}
This implies $x^\eps(t) = x(t) + \eps \xi(t) + O(\eps).$ Since $\xi(t)\geq 0$ we have
  \begin{align}
  u^{\eps,+} = u^\eps(t,x^{\eps,+}(t) ) = (u+\eps v)(t,x(t)+\eps \xi(t)+), \\
  u^{\eps,-} = u^\eps(t,x^{\eps,-}(t) ) = (u+\eps v)(t,x(t)+\eps \xi(t)-).
  \end{align}
Formal  Taylor expansion with respect to $\eps$ shows that 
\begin{align}
 u^{\eps,+}- u^{\eps,-} = \Delta (  u ) + \eps  \Delta ( v+ \xi(t) u_x ) + O(\eps^2), \\
 f(u^{\eps,+}) = f(u^+) + \eps \left( f'(u) ( v+ \xi(t) u_x )\right)^+. 
\end{align}
Therefore, 
\begin{align}
\eps \frac{d}{dt} \xi(t) +O(\eps^2)&= \frac{ \Delta f(u^\eps) }{ \Delta u^\eps }  - \frac{ \Delta f(u) }{ \Delta u } \\&= \frac{ \Delta f(u) +  \eps \Delta ( f'(u) (v+\xi(t) u_x )  \Delta (u) - \Delta f(u) \Delta ( u + \eps ( v + \xi(t) u_x) )  }{ 
	\Delta(u) \left( \Delta(u) + \eps  \Delta ( v + \xi(t) u_x ) \right)
} \\
&= \eps \frac{ \Delta(  f'(u)( v+ \xi(t) u_x)) \Delta(u)- \Delta(f(u)) \Delta(v+\xi(t)u_x)
}{ \Delta(u)^2 }. 
\end{align}
The last equation coincides with equation \eqref{tv-xi2} and this finishes the proof. 
\par 
Lemma \eqref{thm2} implies that for any finite-volume scheme it suffices to include an additional computational step for the shock position $x_i(t)$ as discretization of equation \eqref{RH}. The value of $x_i$ is not necessary to compute the actual solution $u=u(t,x)$ of equation \eqref{claw} 
but required for AD purposes in the sense of tangent vectors \eqref{update}. The details of the implementation are outlined in the following section.

\section{Numerical Method}\label{ad sec}

In this section, we present the numerical method to compute the full tangent required by the new calculus. First, to reduce the abstractness of the presentation we go through the solution procedure referring to the Lax-Friedrichs scheme and a specific problem solved by the Burgers equation. Then, having presented the algebra of the procedure, the implementation combining a flow solver and an AD tool is outlined. We define two methods, i.e., {black-box AD} and {shock AD}.

\subsection{Solution Procedure}

We consider a numerical discretization using finite--volume methods \cite{LeVeque2002aa}. For simplicity we describe  the application in $d=1.$ 
Denote by  $(X_i)_{i \in \mathbb{N}}$  an equidistant spatial grid on $\R$ and  denote by 
$\Delta X= X_{i+1}-X_i.$ The cell boundaries are  $X_{i-\frac12} = X_i - \frac12 \Delta X.$ Then, 
the cell average $\u_i(t)$ on $[X_{i-\frac{1}2}, X_{i+\frac{1}2}]$ at time $t$ for any function $u(t,x)$ is
defined by 
\begin{align}\label{eq::fv_average}
\u_i(t) = \frac{1}{\Delta X}\int_{X_{i-\frac{1}2}}^{X_{i+\frac{1}2} } u(t,x)dx.
\end{align}
 A semi-discretized finite-volume scheme is then given by 
 \begin{align}\label{fv}
 \frac{d}{dt} \u_i(t) = \frac{1}{\Delta X} \left( F_{i +\frac12}(t) - F_{i -\frac12}(t) \right)
 \end{align}
and the initial condition is $\u_{0,i}=\frac{1}{\Delta X}\int_{X_{i-\frac{1}2}}^{X_{i+\frac{1}2} } u_0(x)dx$.  Several choices for the  numerical flux $F_{i +\frac12}(t)$  are known and  we refer to the literature for more details, see e.g. \cite{LeVeque2002aa} and the references therein. The numerical flux $F_{i +\frac12}(t)$ depends on the reconstruction of $\u(t,x)$ at $X_{i+\frac12}$ based on cell averages $\u_j$ for $j\in \mathbb{N}$, see \eqref{recon} below.  In the case
of first--order schemes we have $ j \in \{i-1,i,i+1\}$ and piecewise constant reconstruction of $u(t,x)$ is used.  Furthermore, a suitable time-discretization 
needs to be applied to solve equation \eqref{fv} numerically. As an example for the final fully discrete scheme, we may use the Lax--Friedrichs scheme. The fully discrete form for $\u^n_i=\u_i(t^n)$ reads  for $i \in \mathbb{N}$ and $n \in \mathbb{N}^+:$
\begin{align}\label{scheme}
\u^{n+1}_i = \frac{ \u^n_{i+1}+\u^n_{i-1}}2 + \frac{\Delta t}{2\Delta x} \left(  f(\u^n_{i+1}) - f(\u^n_{i-1})
\right).
\end{align}
The initial data are given by 
\begin{align}\label{IC}
\u^{0}_i = \u_{0,i}.
\end{align}

The time step $\Delta t$ needs to fulfill a Courant-Friedrichs-Lewy (CFL) condition. A dynamical choice is possible and reads for example at  $t=t^n:$
\begin{align}
C_n \Delta t \leq  \Delta X, \mbox{ and }  C_n = \max\limits_{j \in \mathbb{N} }  | f'( \u^n_j )|.
\end{align}
\par 
As discussed before  a numerical approximation to the tangent vector $v(t,x)$ is obtained by AD of the numerical code $G_1$  implementing  equation \eqref{scheme}, i..e, $G_1$ maps $\uad_0:=\left( \u_{0,i} \right)_{i}$ to $\uad_n:=\left(\u^n_i\right)_i.$ 
In order to illustrate the AD we give the respective AD of the Lax--Friedrichs scheme \eqref{scheme}
as 
\begin{align}\label{uwe_ad1}
\dot{\u}^{n+1}_i=\frac{\dot{\u}^n_{i+1}+\dot{\u}^n_{i-1}}2 + \frac{\Delta t}{2\Delta X} \left( \frac{d}{d \u} f(\u^n_{i+1}) \cdot \dot{\u}^n_{i+1} - \frac{d}{d \u} f(\u^n_{i-1}) \cdot \dot{\u}^n_{i-1}
\right),
\end{align} 
where $\dot{\uad}_n=\left(  \dot{\u}^n_i\right)_i$. Clearly, we expect $\dot{\u}^n_i$ to be an approximation 
of the cell average  $$ \dot{\u}^n_i \approx \frac{1}{\Delta X} \int_{X_{i-\frac{1}2}}^{X_{i+\frac12}}  v(t^n,x) dx$$ where $v$ is the solution to equation \eqref{tv-v}. Due to the previous theoretical discussion  the   knowledge of $\dot{\u}^n_i$ to be not sufficient to characterize the full tangent vector. 
\par 
In order to obtain a  numerical approximation to the full  tangent vector $(v,\xi)$ the finite--volume scheme \eqref{fv} is augmented by $i=1,\dots, N(u_0)$ additional equations where $N(u_0)$ denotes the number of discontinuities in $u_0(\cdot)$, see  Lemma \eqref{thm}. The initial  position of the discontinuities  are denoted by $x_i(0)=x_{i,0}$ for $i=1,\dots, N(u_0).$  In the following we discuss a  numerical discretization of equation 
\begin{align}\label{rh temp}
\frac{d}{dt} x_i(t) = \frac{ \Delta f(u) }{\Delta u}, \; x_i(0)=x_{i,0}, \; i=1,\dots,N(u_0)
\end{align}
such that  AD leads to a consistent approximation for the evolution of $\xi_i(t)$ where $\xi(t)$ is obtained by equation \eqref{tv-xi2}. 
On the continuous level  the linearization of  equation \eqref{rh temp} is consistent with equation \eqref{tv-xi2}. However, a straight--forward
numerical discretization of equation \eqref{rh temp} e.g. given by an explicit Euler scheme leads to the numerical approximation 
$\x^n_i $ of $x_i(t^n)$ by  
\begin{align}\label{wrong shock discr}
\x^{n+1}_i = \x^n_i + \Delta t \left( 
\frac{f(\u(t^n,\x_i^{n}+))-f(\u(t^n,\x_i^{n}-))}{\u(t^n,\x_i^{n}+)-\u(t^n,\x_i^{n}-)}
\right), \;  x^0_i = x_{i,0}.
\end{align}
Here, we denote by $\u(t,\cdot)$ the piecewise constant reconstruction on $\R$ based on the cell averages $(\u_i^n)_{i}^n$ for $i \in \mathbb{N}$ and $ n\geq 0$:
\begin{align}\label{recon}
\u(t,x) = \sum\limits_{ n\geq 0} \sum\limits_{i\in\mathbb{N}} \u_i^n \; \chi_{[t^n,t^{n+1}]}(t)  \chi_{ [X_{i-\frac12}, X_{i+\frac12} ]} (x).
\end{align}
The {black-box} application of tangent AD would then yield the  numerical approximation $\dot{\x}_i^n$
 to the tangent $\xi_i(t)$ at $t=t^n$ for each $i$ as  
\begin{equation} \label{eqn:wrong}
\dot{\x}^{n+1}_i =\dot{\x}^n_i+\Delta t \cdot \frac{d}{d \x} 
\frac{f(\u(t^n,\x_i^{n}+))-f(\u(t,\x_i^{n}-))}{\u(t,\x_i^{n}+)-\u(t,\x_i^{n}-)} \cdot \dot{\x}^n_i,\; \dot{\x}^0_i = \xi_{i,0}.
\end{equation}
However, this approach yields potentially wrong approximations to the tangent vector, i.e.,  
\begin{align*}
\dot{\x}^n \not = \xi(t^n) +O(\Delta t)^2
\end{align*}
 The reason is two-fold: any numerical finite--volume scheme of the type \eqref{fv} introduces artificial viscosity. In the case of a first-order scheme this introduces
an error of  $O(\Delta X^2).$ This prevents a sufficiently sharp resolution of the shock. 
Further, equation \eqref{rh temp} is numerically unstable if $ \u(t^n,\x^n_i + ) \approx   \u(t^n,\x^n_i- ).$
 This might occur in the numerical scheme even if there is a shock located at $\x^n_i$ in the case when
 this  location  is not sufficiently sharply resolved.  
 Second, for a discretization of $\Delta \u$ we require a reconstruction of $(t,x) \to \u(t,x)$ and an evaluation at $\x_i^n\pm.$
 For a first--order numerical scheme, $\u$ is reconstructed piecewise constant as given by equation \eqref{recon}. Hence, 
 black--box AD applied to equation \eqref{wrong shock discr} will not be able to recover numerical approximation to the  terms $\partial_x u^\pm$
 of equation \eqref{tv-xi2}. Hence, solely for the purpose of applying  AD, we  propose  one--sided piecewise linear reconstruction of $x\to u(t,x)$ in the vicinity of $x_i(t).$
  \par 
 In order to address both points we therefore implement a numerical approximation to equation \eqref{RH} depending on two parameters $C>0, \alpha \geq 1$ as follows. We assume $(\u^n_i)_i$ are given by any  finite volume scheme \eqref{fv} as for Example \eqref{scheme}.  
 Initial  positions for shocks  are given by $x_i(0)$ for $i=1,\dots, N(u_0)$. We proceed using the following steps:
 \begin{itemize}
 	\item Use the piecewise constant reconstruction $\u$ \eqref{recon} and approximate the shock position $x_i(t^n)$ by $\x_i^n$ given by  
 	\begin{align}\label{shock position num}
 	\x^{n+1}_i = \x^n_i + \Delta t   \frac{d}{d\u} f( \u(t^n,\x^n_i )) , \; \x^0_i= x_i(0). 
 	\end{align}

\item Additionally, consider  a  piecewise linear reconstruction $\v(t,x)$ of $u(t,x)$, i.e., 
\begin{align}\label{recon linear}
\v(t,x) := \sum\limits_{ n\geq  0} \sum\limits_{i \in \mathbb{N} } \left(  \u_i^n + (x-X_i) \us^n_i \right) \chi_{[t^n,t^{n+1}]}(t) 
 \chi_{ [X_{i-\frac12}, X_{i+\frac12} ]} (x). 
\end{align}
Several possibilities for approximation $\us^n_i$ exist. For example, non--oscillatory reconstructions can be used \cite{LeVeque2002aa}. 
 Since we have not seen  any major improvement using e.g. a reconstruction using a minmod limiter compared with the following one--sided differences:
\begin{align}
({\bf U_x^+})^n_i = \frac{1}{\Delta X} \left( \u^n_{i+1} - \u^n_i\right), \; ({\bf U_x^-})^n_i =  \frac{1}{\Delta X}  \left(\u^n_i -  \u^n_{i-1} \right).
\end{align}
Note that those slopes are {\em not} used to propagate the solution $\u^n_i$. 
The reconstruction $\v$ is used only as an auxiliary variable used to allow for a numerical approximation 
to the Rankine--Hugenoit condition such that AD is applicable. Hence, the next step consists of replacing equation \eqref{wrong shock discr}. 

\item Consider the approximation of the Rankine--Hugenoit condition replacing \eqref{wrong shock discr}. 
 \begin{align}\label{second shock position num}
{\bf y}^{n+1}_i  := \x^n_i + \Delta t \; \frac{  f(\v(t^n, \x_i^n + C (\Delta X)^\alpha )  )    -  f(\v(t^n, \x_i^n - C (\Delta X)^\alpha ))      }{  
	\v(t^n, \x_i^n + C (\Delta X)^\alpha )  - \v(t^n, \x_i^n - C (\Delta X)^\alpha )   }  
\end{align}
Again, ${\bf y}^{n+1}_i$ is an auxiliary  variable used to apply the AD framework.  Hence, we approximate $\xi_i(t^n)$ by $\dot{{\bf y}}^n_i$, see equation \eqref{SRH}.  Note that this implies we define $\dot{\x}^{n+1}_i$ as the AD tangent by the differentiation of equation \eqref{second shock position num}.  
 \end{itemize}
Hence, for applying AD we separate the evolution of the shock position \eqref{shock position num} and  computation
of a suitable tangent using AD  according to equation \eqref{second shock position num}.   Summarizing, 
the proposed procedure leads to the following set of equations for $n \geq0$ and $i=1,\dots,N(u_0)$ 
\begin{align}\label{uwe_ad2_1}
\x^{n+1}_i &=\x^n_i+\Delta t \cdot \frac{d}{d \u} f(\u(t^n,\x_i^{n})),\ \x_i(0)=x_{i,0}\\ 
\label{SRH}
\dot{\x}^{n+1}_i :=\dot{{\bf y}}^{n+1}_i&=\dot{\x}^n_i+\Delta t \cdot \frac{d}{d \x} 
\frac{f(\v(t^n,\x_i^{n,+}))-f(\v(t^n,\x_i^{n,-}))}{\v(t^n,\x_i^{n,+})-\v(t^n,\x_i^{n,-})} \cdot \dot{\x}^n_i,\ \dot{\x}^0_i=\xi_{i,0}\\
\x^{n,\pm}_i  &:=  \x_i^n \pm \srw \label{eq::shock_delta}
\end{align}
with the half width of the numerical approximation of the shock $\srw = C (\Delta X)^\alpha$.

Some remarks are in order.  We expect that $\dot{\u}^{n+1}_i  \approx v(t^n,x_i)$ and $\dot{\x}^n_i\approx \xi_i(t^n)$ where $(v,\xi)$ are the tangent vectors introduced above.  Here, $\dot{\x}^n_i$ is the AD formulation  applied to equation \eqref{second shock position num} leading formally to  equation \eqref{SRH}.   The particular discretization \eqref{shock position num}  is an approximation to the continuous formulation \eqref{RH}. It is proposed to provide a remedy to the drawbacks of a  straightforward discretization of equation \eqref{RH}. 
 Since the artificial diffusion is of order $\Delta X^2$ a possible choice
for $C$ and $\alpha$ would be $C=1$ and $\alpha=2.$ Clearly, if  $\x^n_i \pm \srw$ is outside
the numerical approximation of the shock, the previous formula does not provide an approximation 
to the true propagation speed of the shock. It is also clear that the particular choice of $C$ and $\alpha$ depend 
on the underlying finite-volume  scheme and  the applied reconstruction procedure. So far, 
we can not provide a general formula for choosing $C$ and $\alpha$. 

\subsection{Implementation}\label{susec::implementation}

So far, the details of the numerical procedure were exemplified for the Lax-Friedrichs scheme \eqref{scheme}. Keeping the computation of tangents to real-world numerical simulations in mind, the essential routines of a state-of-the-art CFD code, i.e., the Zonal Flow Solver (ZFS), are adapted to the proposed procedure. ZFS is a highly efficient multi-physics simulation framework and developed by the numerical group of the Institute of Aerodynamics and Chair of Fluid Mechanics (AIA), RWTH Aachen University.
Flows with shocks were computed for several applications, i.e., a transonic airfoil in \cite{schneiders13}, a cone in supersonic flow in \cite{schilden15}, a blunt stagnation point probe in supersonic flow in \cite{schilden17,schilden18}, and a reentry capsule in supersonic flow in \cite{schilden20}.

Here, the computational domain of the numerical simulation is discretized by an unstructured Cartesian grid and the governing equations are integrated using a finite-volume method~\cite{lintermann14,hartmann08}. For the spatial discretization, an advection upstream splitting method (AUSM) is used. The cell center gradients are computed using a second-order accurate least-squares reconstruction scheme~\cite{schneiders16}. Shock capturing is achieved by adding additional numerical dissipation at the shock position using a slope limiter~\cite{venkatakrishnan93}. The temporal integration is based on a 5-stage second-order accurate Runge-Kutta scheme.
Three steps are required to obtain the full tangent $(v,\xi)$ from the numerical simulation.
First, the CFD code is modified to solve the additional equation for the shock location \eqref{uwe_ad2_1}.  Second, the modified method is enhanced by AD to yield the first component of the tangent $v$, i.e., the numerical discretization of equation \eqref{tv-v}. Finally, we customize the numerical computation of the Rankine--Hugenoit condition  \eqref{second shock position num} such that AD provides the second component of the tangent $\xi_i$. 

The computation of the full tangent, i.e., an equation equivalent to \eqref{uwe_ad1} and equation \eqref{SRH}, is based on the AD software tool dco/c++\footnote{dco/c++ is developed by the Numerical Algorithms Group Ltd. in collaboration with the STCE group at RWTH Aachen University; see also \url{https://www.nag.co.uk/content/algorithmic-differentiation-software}.}\cite{AIB-2016-08}. 
It relies on function and operator overloading in combination with 
extensive C++ template metaprogramming to yield highly efficient tangent 
(and adjoint) code of arbitrary order. Program variables with non-vanishing derivatives are re--declared as {\em active}. Corresponding tangent and adjoint types are provided by dco/c++. The set of {\em elemental} functions including all built-in arithmetic operators and intrinsic functions is 
overloaded for the active (for example, first-order tangent) data type. 
dco/c++ has been applied successfully to numerous practically relevant applications in
Computational Science, Engineering, and Finance; see, for example, \cite{Towara2018Sam,Naumann2018AAD,Maybank2020MfB}.

Real-world numerical simulations subject to parameter sensitivity 
analysis, nonlinear optimization or optimal control often require 
another approach than {black-box AD}. Typical reasons include infeasible persistent memory requirement in adjoint mode, calls to binary third-party library functions as well as local nondifferentiability -- the later
 is also part of the problem tackled in this paper. 
Solutions with dco/c++ rely on the extension of the set of elemental 
functions with solutions for the respective subproblems. 
The application of AD is locally replaced by 
a call to a specifically designed method. In the given context the 
{black-box} application of AD to the evolution of the shock location in equation~\eqref{eqn:wrong} is replaced by a custom elemental function 
implementing equations \eqref{uwe_ad1} and 
\eqref{SRH}. dco/c++ treats it similar to any 
other built-in function. Abstraction is lifted 
to the level necessary for dealing with the discontinuity due to the 
shock in a numerically consistent way. This approach enables correct 
approximation of tangents as outlined above and is termed {shock AD}. 
A comprehensive discussion of the software engineering aspects is beyond the scope of this contribution.

Algorithmic adjoint parameter 
sensitivities of the shock location within ${\bf x}$ follow seamlessly. 
Their
implementation with dco/c++ uses a corresponding custom adjoint 
elemental function. The extension to algorithmic adjoints of (objectives defined over the final) state ${\bf U}$ turns out to be less straightforward  in general due to nonlinearity $f.$ This is the subject of ongoing research.

Summarizing \S~\ref{ad sec}, the complete scheme is given by the fully discrete finite-volume scheme as for example given by \eqref{scheme}, 
the propagation of the shock position \eqref{shock position num} and based on this position the update \eqref{second shock position num}. The AD is then given by equations \eqref{uwe_ad1} and \eqref{SRH}, respectively. 
Next, we present the numerical results applying the new calculus to problems for the Burgers and Euler equations.
\section{Computational Results} \label{sec:cr}
In this section, the results applying the theoretical calculus are presented. First, we continue with the example of the Burgers equation that supported the presentation of the calculus in \S~\ref{sec::calculus}. We compare the methods denoted by {black-box AD} and {shock AD}. Finally, we refer to an example of the Euler equations.
\subsection{Tangent Vectors for Burgers Equation}\label{susec::cr::burgers}
\begin{figure}
\centering
\subfloat[\label{sufig::burgers::solution_u}]{
\includegraphics[scale=\figscale]{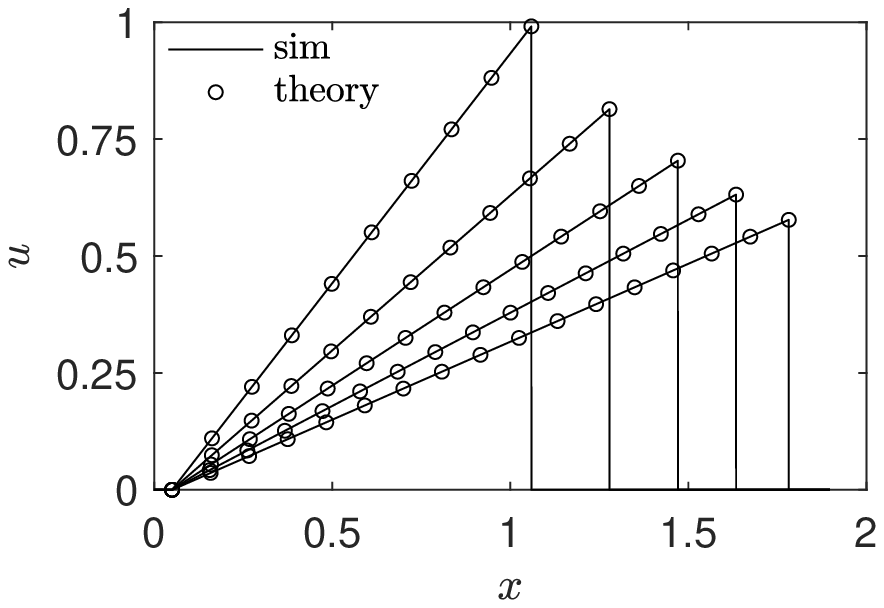}
}
\qquad
\subfloat[\label{sufig::burgers::solution_xs}]{
\includegraphics[scale=\figscale]{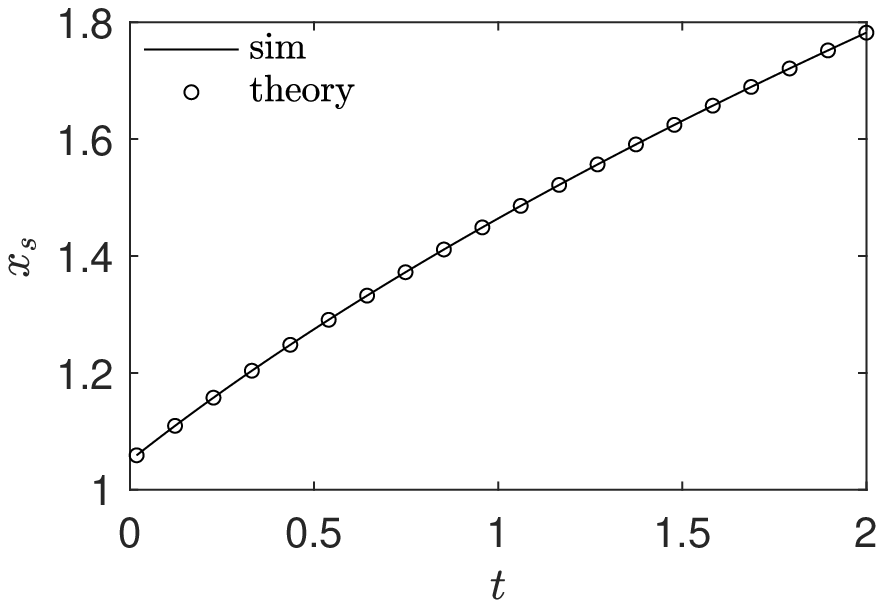}
}
\caption{Comparison of numerical and theoretical solutions for Example \ref{ex::burgers::ramp}: (a) Burgers solution at $t \approx (0.018, 0.51, 1, 1.5, 2)$; (b) shock location.}\label{fig::burgers::solution}
\end{figure}

\begin{figure}
\centering
\subfloat[\label{sufig::burgers::derivative_u}]{
\includegraphics[scale=\figscale]{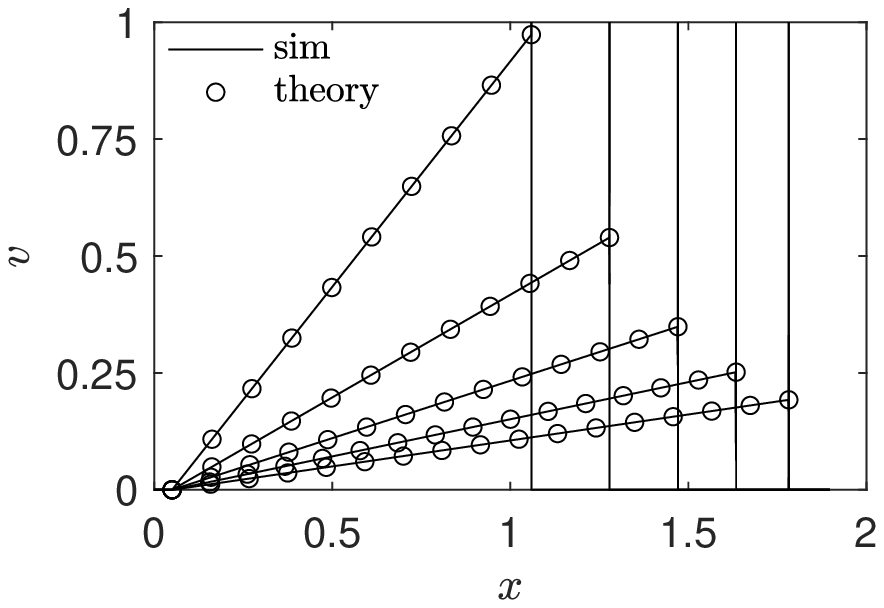}
}
\quad
\subfloat[\label{sufig::burgers::derivative_xs}]{
\includegraphics[scale=\figscale]{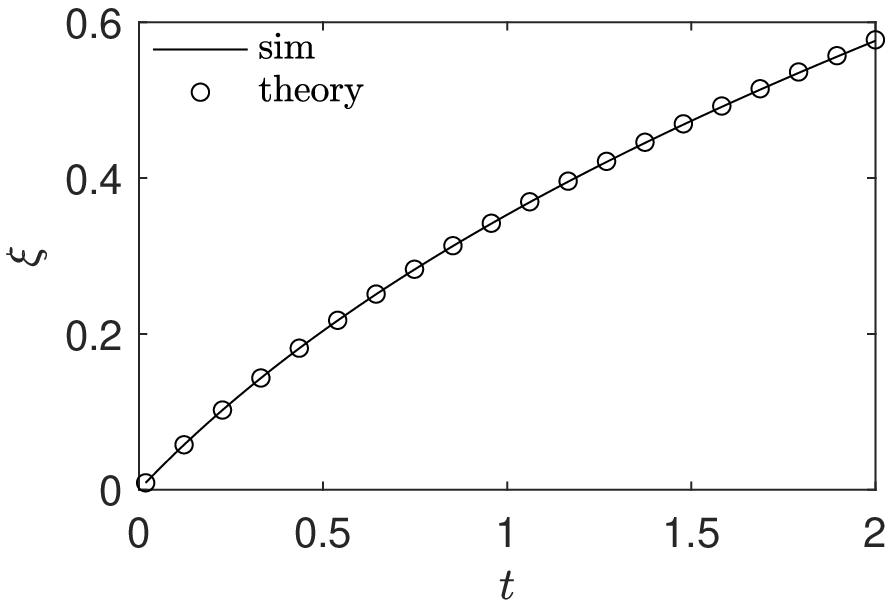}
}
\caption{Comparison of the tangent of the numerical and theoretical solutions for Example \ref{ex::burgers::ramp}: (a) Tangent to Burgers solution at $t \approx (0.018, 0.51, 1, 1.5, 2)$; (b) tangent to shock location.}\label{fig::burgers::derivative}
\end{figure}

Example \ref{ex::burgers::ramp} of the Burgers equation \eqref{eq::burgers} is solved by the numerical method presented in \S~\ref{susec::implementation}. Note that the initial condition \eqref{eq::ramp::ic} is shifted by $0.05$ to greater $x$.
The example is computed on nine equidistant grids. From grid to grid with increasing No.\ in Table \ref{tab::cases::burgers}, the cell width doubles. Depending on the grid, the computational domain extends over $0 \leq x \lesssim 1.9$. To yield a final solution time of $t_{final} = 2$ using a CFL number of $\approx 0.63$, the constant time steps $\Delta t$ of Table \ref{tab::cases::burgers} are chosen.
\begin{table}
\centering
\resizebox{\textwidth}{!}{
\begin{tabular}{lccccccccc}
No.       &   1  &  2  &  3  &  4  &  5  &  6  &  7  &  8  &  9\\
\hline
$\Delta X$& $5.75\mbox{e}^{-5\strut}$&$1.15\mbox{e}^{-4}$&$2.3\mbox{e}^{-4}$&$4.6\mbox{e}^{-4}$&$9.2\mbox{e}^{-4}$&$1.84\mbox{e}^{-3}$&$3.68\mbox{e}^{-3}$&$7.36\mbox{e}^{-3}$&$1.472\mbox{e}^{-2}$\\
\hline 
$\Delta t$ & $3.64\mbox{e}^{-5\strut}$  & $7.27\mbox{e}^{-5}$  & $1.45\mbox{e}^{-4}$  &  $2.9\mbox{e}^{-4}$  & $5.88\mbox{e}^{-4}$  & $1.18\mbox{e}^{-3}$  & $2.35\mbox{e}^{-3}$  & $4.7\mbox{e}^{-3}$  & $9.52\mbox{e}^{-3}$ \\
\hline 
\end{tabular}}
\caption{Parameters of the simulations of Example \ref{ex::burgers::ramp} of the Burgers equations \eqref{eq::burgers}: No.\ is the reference to the simulation, $\Delta X$ is the cell width, and $\Delta t$ is the time step with three significant digits. Common parameters are the CFL number $cfl\approx 0.63$, and $C=5$, and $\alpha=1$ in \eqref{eq::shock_delta}.}\label{tab::cases::burgers}
\end{table}
First, we validate the results of simulation No.\ 1 using  analytical results of Examples \ref{ex::burgers::ramp} and \ref{ex::burgers::ramp::solution}, i.e., we compare Burgers solution $u$ in \fig{\ref{sufig::burgers::solution_u}} to \eqref{eq::burgers::solution_u} for $\eps=0$, the shock position $x_s$ in \fig{\ref{sufig::burgers::solution_xs}} to $x(t)$ in \eqref{eq::burgers::shock_location}, the tangent to the Burgers solution $v$ in \fig{\ref{sufig::burgers::derivative_u}} to \eqref{v1}, and the tangent to the shock location $\xi$ in \fig{\ref{sufig::burgers::derivative_xs}} to \eqref{xi1}.
The agreement of simulation and theory is evident. Relevant deviations occur in vicinity of the shock defined by \eqref{eq::shock_delta}. First, the shock is smeared over several cells. Second, this continuous numerical representation of the discontinuity leads to tremendous deviations of the tangent vector $(v,\xi)$ to the solution in \fig{\ref{sufig::burgers::derivative_u}} that need to be considered in the computation of the calculus.

In terms of the calculus for the sensitivity, the data in \fig{\ref{fig::burgers::solution}} denotes the function evaluation, whereas the data in \fig{\ref{fig::burgers::derivative}} is the tangent $(v,\xi).$ Using the function value and the  tangent vector with respect to $\eps$, a function value $\tilde{u}$ for $\eps \neq 0$ can be approximated by the methods defined in \S~\ref{susec::implementation}, i.e., {black-box AD} and {shock AD}. The approximation is also called tangential shift since the solution is shifted in the direction of the tangent. Considering the shock location, we refer to a tangential displacement. The characteristic function bridges the gap from the location of the shock in the simulation to the approximated position for $\eps \neq 0$ and thus, the shock is displaced. 

\begin{figure}
\centering
\includegraphics[scale=0.85]{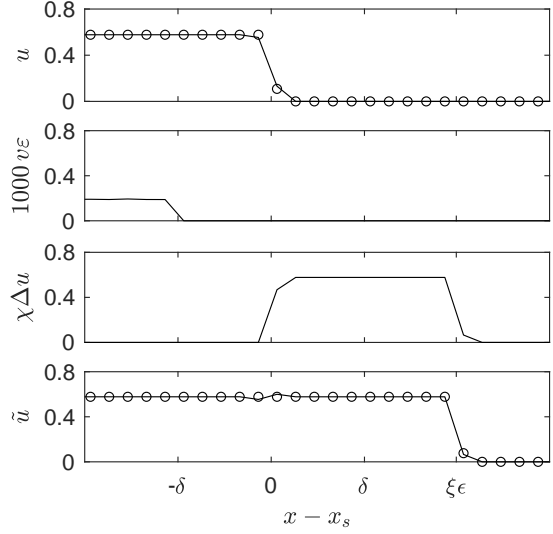}
\caption{Detailed view of the components of the calculus \eqref{eq::calculus} in the vicinity of the shock at $x_s$: Solution of the Burgers equation $u$, tangential shift $v\eps$ for $|x-x_s| > \srw$, the tangential displacement of the shock $\chi\Delta u$ following \eqref{eq::fv_average}, and the final approximation $\tilde{u}$. The symbols denote the analytic solution \eqref{eq::burgers::solution_u} in compliance with \eqref{eq::fv_average}. For this illustration the perturbation is $\eps\approx 1\mathrm{e}^{-3}$.}\label{fig::burgers::calculus}
\end{figure}

The application of {shock AD} employing the calculus \eqref{eq::calculus} is not straightforward and requires two more steps. To remove the erroneous tangent to the Burgers solution in the vicinity of the shock, we omit the second term on the right-hand side of \eqref{eq::calculus} in the region specified by \eqref{eq::shock_delta}. For the third term, the characteristic function $\chi$ has to be evaluated in the discretized computational domain. We follow \eqref{eq::fv_average} for $\chi$ to generate an approximation $\tilde{u}$ that is continuous in $\eps$. The terms of the calculus \eqref{eq::calculus} in the vicinity of the shock at $x_s$ are shown in \fig{\ref{fig::burgers::calculus}}. Note that the coordinate is relative to the shock at $x_s$.
The symbols denote the analytic solution.
Again, computing the analytic solution, we use cell averages \eqref{eq::fv_average}.
The shock is in good agreement at both locations, i.e., at $x_s$ and $x_s+\xi\eps$. The contribution of the erroneous tangent $v$ to the Burgers solution in $-\srw < x-x_s < \srw$ is omitted.

\begin{figure}
\centering
\subfloat[\label{sufig::burgers::convergence::perturbation}]{
\includegraphics[scale=\figscale]{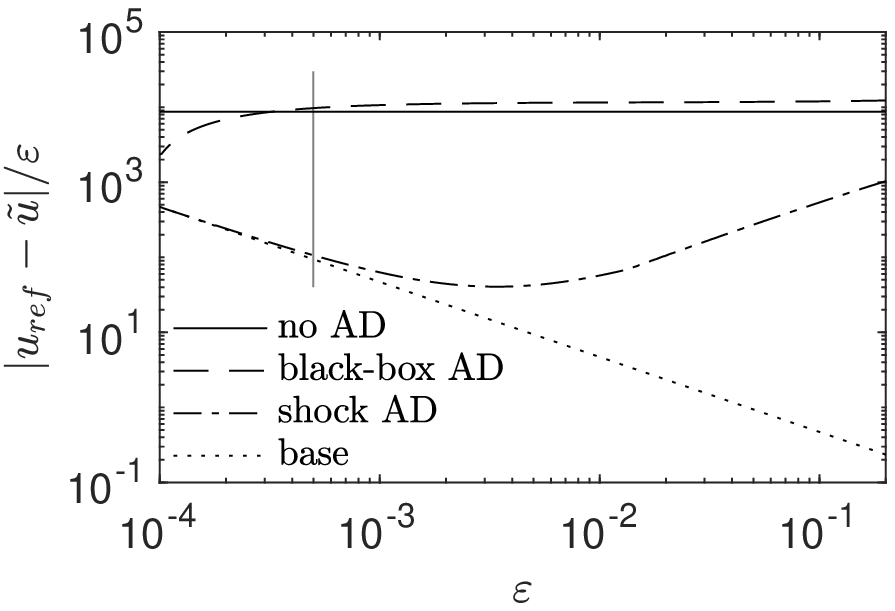}
}
\quad
\subfloat[\label{sufig::burgers::convergence::grid}]{
\includegraphics[scale=\figscale]{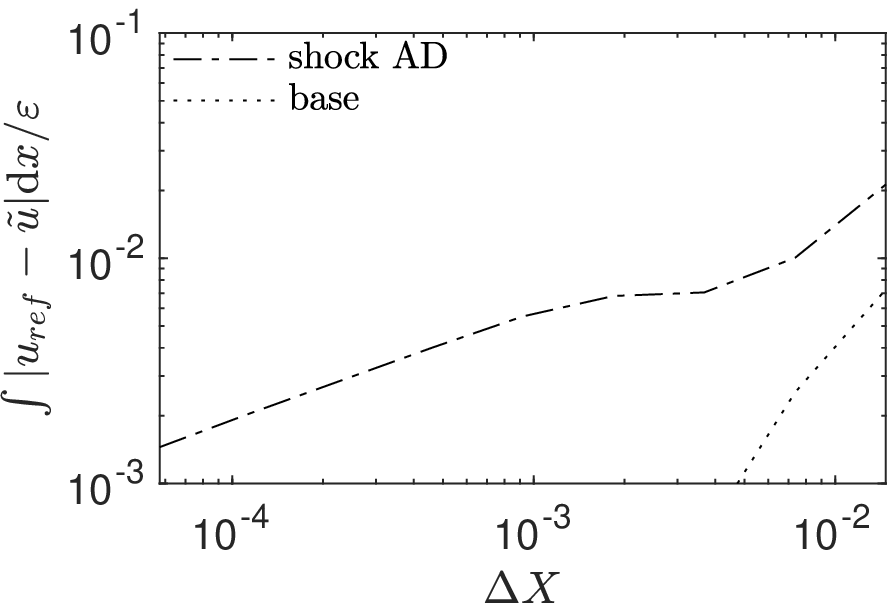}
}
\caption{Convergence of the tangential shift for Example \ref{ex::burgers::ramp}: (a) Convergence with respect to the perturbation $\eps$;
the tangential shift is compared to the analytic solution applying
``no AD'' no tangential shift, 
``block-box AD'' tangential shift of solution $u$,
``shock AD'' tangential shift using the full tangent $(v,\xi)$,
``base'' reference error of the simulation divided by $\eps$;
(b) Convergence with respect to the grid spacing $\Delta x$; perturbation is $\eps_{max}$;
number of cells $n_{min} = 129 \approx n_{max}/256$ to $n_{max}=33000$;
``var'' denotes the error of the tangential shift;
``base'' denotes the error of the solution without tangential shift.
}\label{fig::burgers::convergence}
\end{figure}

To validate the calculus as an adequate measure for the tangent to the solution with respect to $\eps$, the convergence of the error of the approximation based on the $L_1$ norm and the analytic solution \eqref{eq::burgers::solution_u} is analyzed in \fig{\ref{fig::burgers::convergence}}.  
First in \fig{\ref{sufig::burgers::convergence::perturbation}}, the convergence of the error of the approximation with respect to $\eps$ is shown using case No.\ 1. The maximum perturbation $\eps_{max} = 0.2$ shifts the shock approximately to the end of the computational domain.
The minimum perturbation of $\eps_{min} = 1\mathrm{e}^{-4}$ is defined by $\eps_{min} = \Delta X/\xi$, i.e., the tangential displacement of the shock is equal to the grid spacing. Furthermore, $\eps_\dagger = \srw/\xi$ is denoted by a vertical line. For $\eps < \eps_\dagger$, the displaced shock is in the region of the grid that contains the numerical representation of the shock.

The deviation of the simulation ``base'' from the analytic solution for $\eps=0$ is constant and thus, the data increases with decreasing $\eps$. This is the reference error caused by the discretization. Using no means of sensitivity to approximate $u_{ref}$ leads to $\tilde{u} = u$. The error with respect to $\eps$ is nearly constant. The black-box AD method shows a similar behavior for $\eps_{\dagger} \leq \eps \leq \eps_{max}$. The error of the proposed shock AD method converges toward the reference error ``base.'' Note that the approximation is based on the function evaluation in \fig{\ref{fig::burgers::solution}} and the full tangent in \fig{\ref{fig::burgers::derivative}}. For $\eps \rightarrow 0$, an error smaller than ``base'' can not be expected.
This constitutes the main result of this paper. The {black-box AD} application yields meaningless sensitivities, whereas the full tangent of the {shock AD} method provides the correct sensitivities to achieve the desired error convergence in \fig{\ref{sufig::burgers::convergence::perturbation}}. 

The decrease of the black-box AD data for $\eps \leq \eps_{\dagger}$, is artificial. Numerical dissipation causes a continuous shock in the region $x_s - \srw < x < x_s + \srw$ and if the shock displacement is small, the derivatives in this region give a meaningful approximation. Note that this behavior is based on numerical dissipation and not on the theory of hyperbolic partial differential equations with shocks.

The analysis in \fig{\ref{sufig::burgers::convergence::perturbation}} shows that the error with respect to $\eps$ of the shock AD method converges toward the reference error ``base.'' For $\eps \rightarrow 0$, convergence of the error toward zero can not be shown when the method is compared to the analytic solution. This effect is caused by the error due to the discretization and not by the calculus that defines the {shock AD} method. By showing grid convergence of the new method, this flaw can be remedied. In \fig{\ref{sufig::burgers::convergence::grid}}, the error of the shock AD method and the reference error ``base'' are shown for an approximation using a constant perturbation of $\eps = \eps_{max}$. Both vanish with increasing grid refinement and thus, the minimum of the shock AD error in \fig{\ref{sufig::burgers::convergence::perturbation}} reduces and vanishes for $\Delta x \rightarrow 0$. That is, the implemented {shock AD} method yields sensitivities consistent with the theory in \S~\ref{sec::calculus} based on \eqref{eq::calculus}.
Next, we present an example for the Euler equations.
\subsection{Tangent Vectors For Euler Equations}\label{susec::cr::euler}
Before the results are analyzed, the Euler equations and the computational setup are presented. The problem statement is complemented by an example to exemplify the involved physics.
\subsubsection{Problem Statement}
The Euler equations are a hyperbolic system of partial differential equations which describe the conservation of mass, momentum, and energy of an inviscid fluid and read
\begin{equation}\label{eq::euler}
\int\limits_V \frac{\partial \mathbf{Q}_c}{\partial t} \mathrm{d} V + \oint\limits_A \overline{\mathbf{H}} \cdot \mathbf{n} \mathrm{d} A = \mathbf{0}\,.
\end{equation}
The quantity $\mathbf{Q}_c = \left[ \rho, \rho\mathbf{u}^T, \rho E \right]^T$ is the vector of the conservative variables with the density $\rho$, velocity vector $\mathbf{u}$, and the total specific energy $E = e + \mathbf{u}^2/2$ containing the specific energy $e$.
Alternatively, the flow can be described by the primitive variables $\mathbf{Q}_p = \left[ \rho, \mathbf{u}^T, p \right]^T$, where $p$ denotes the static pressure.
The flux vector $\overline{\mathbf{H}}$ is
\begin{equation}\label{eq::euler::flux}
\overline{\mathbf{H}} = \ar{c}{\rho \mathbf{u}^T\\\rho\mathbf{u}\mathbf{u}+\overline{\mathbf{I}}p\\\mathbf{u}^T(\rho E + p)}\,.
\end{equation}
The system of equations is closed by the equation of state for an ideal gas
\begin{equation}
e = \frac{p}{\rho\left( \gamma - 1 \right)}
\end{equation}
with the ratio of specific heats $\gamma$. 

Flows governed by the Euler equations may contain shocks, if the Mach number, i.e., the ratio of the flow velocity $u$ and speed of sound $a$, $\Mach = u/a$, is $\Mach > 1$. Then, all characteristics point forward in the direction of the flow and shocks provide the only mechanism to propagate information, e.g., the presence of a body, upstream. This mechanism causes entropy production.

\begin{figure}
\centering
\begin{tikzpicture}[scale=1.5,decoration=snake]
\node[anchor=south west] at (0,0){\includegraphics[scale=1.0]{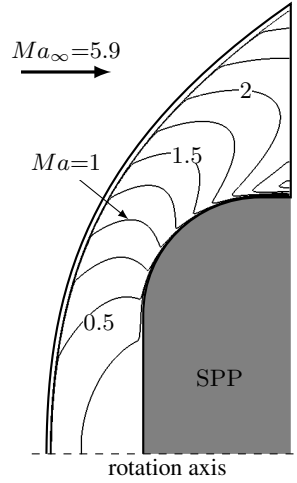}};
\draw[-latex,very thick] (0,3.5) -- +(0.8,0) node[midway,above]{\footnotesize$\Mach_\infty{=}5.9$};
\node[fill=white] at (0.4,2.7){\footnotesize$\Mach{=}1$};
\node at (1.3,0){\footnotesize rotation axis};
\end{tikzpicture}
\caption{Lines of constant Mach number between a detached shock wave and a blunt body, e.g., defined by a stagnation point probe (SPP)~\cite[\textcopyright 2017 Cambridge University Press, Reprinted with permission]{schilden17,schilden18}.}\label{fig::spp}
\end{figure}

An example is shown in \fig{\ref{fig::spp}}. A blunt body, e.g., defined by a stagnation point probe (SPP), is exposed to supersonic flow in a wind tunnel to measure stagnation pressure fluctuations caused by flow perturbations~\cite{schilden17,schilden18}.
The flow is characterized by lines of constant Mach number and a detached shock wave. 
The outer black line illustrates the boundary of the computational domain.
In the freestream field at $\Mach_\infty$, the lines of constant Mach number collapse and denote the location of the shock wave. At the rotation axis, it is a normal shock.
There, the pressure rise of the fluid is most intense and the velocity immediately downstream of the shock is minimum. Transmitting the shock wave, the freestream is not deflected. Thus, it can be modeled by the one-dimensional Euler equations.

Here, we compute the full tangent to the propagation of a normal shock defined by a Riemann problem of the Euler equations.
The parameters of the flow are the initial location of the shock $x_{s,0}$, the Mach number $\Mach$, and the shock speed $S$. The left side of the Riemann problem is defined by the Mach number. The right side follows from the shock speed given that only a single shock occurs. The quantities density $\varrho$, pressure $p$, temperature $T=\gamma p/\varrho$, and velocity $u$ are non-dimensional. The dimensional reference is defined by the temperature, speed of sound, and density at rest. Then, the temperature on the left side of the Riemann problem is
\begin{align}
T_l = \left(1+\frac12 (\gamma-1) \Mach^2\right)^{-1}\,.
\end{align}
The velocity $u_l$, static pressure $p_l$, and density $\varrho_l$ are
\begin{align}\label{eq::euler::left}
\ar{c}{u_l\\p_l\\\varrho_l} = \ar{c}{\Mach \sqrt{T}\\T^{\gamma/(\gammaMinusOne)}/\gamma\\T^{1/(\gammaMinusOne)}}.
\end{align}
The ratio of specific heats for air is $\gamma=1.4$. The left and right states are coupled by the Rankine-Hugoniot conditions~\cite{toro99}.
The solution can be given by ratios of left and right variables. In these ratios
\begin{align}
  \frac{p_r}{p_l} &= 1+\frac{2\gamma}{\gamma+1}\left(\widetilde{\Mach}^2 - 1 \right)\label{eq::euler::pressureratio}\\
  \frac{\varrho_r}{\varrho_l} &= \frac{(\gamma+1)\widetilde{\Mach}^2}{(\gamma-1)\widetilde{\Mach}^2 + 2}\\
  \frac{u_r}{u_l} &= \frac{\varrho_l}{\varrho_r}\,,
\end{align}
the Mach number relative to the moving shock, i.e., $\widetilde{\Mach} = \frac{u_l-S}{a_l}$, is introduced.
The following two equations are essential for the solution algorithm. First, the characteristic defined by the speed
\begin{align}\label{eq::euler::fp}
  \left.\frac{\mathrm{d}x}{\mathrm{d}t}\right|_{sa} = u - a
\end{align}
describes the propagation of slow acoustic waves (sa). For orthogonal shocks, these waves propagate in the opposite direction of the flow and always run into the shock. The shock speed can be computed from, e.g., (\ref{eq::euler::pressureratio}) and is
\begin{align}\label{eq::euler::rh}
  S = u_l - a_l\sqrt{\frac{\gamma+1}{2\gamma}\frac{p_r}{p_l}+\frac{\gamma-1}{2\gamma}}\,.
\end{align}
Analogous to the Burgers equation, \eqref{eq::euler::fp} is used to integrate the shock location and \eqref{eq::euler::rh} is differentiated by the {shock AD} method.

A single simulation result is presented. The computational domain is defined by $0 \leq x \leq 210$ and the cell width is $\Delta X = 0.01$. The initial shock is located at $x_{1,0}=x_{s,0} = 5$. The Mach number is $\Mach = 5.3452$ and the shock speed is $S=0.1$. The final time of the simulation is $t_{final}=1000$ and the CFL number is $cfl = 0.82$. For the {shock AD} method, $C=20$ and $\alpha=1$ are specified.
Tangent vectors with respect to the shock speed are evaluated such that the quantity $\eps$ is a perturbation of $S$. That is, $x_s^\eps=x_{s,0}+(S+\eps)t$ and $\xi = t$.
\subsubsection{Results}
\begin{figure}
\centering
\includegraphics[scale=0.85]{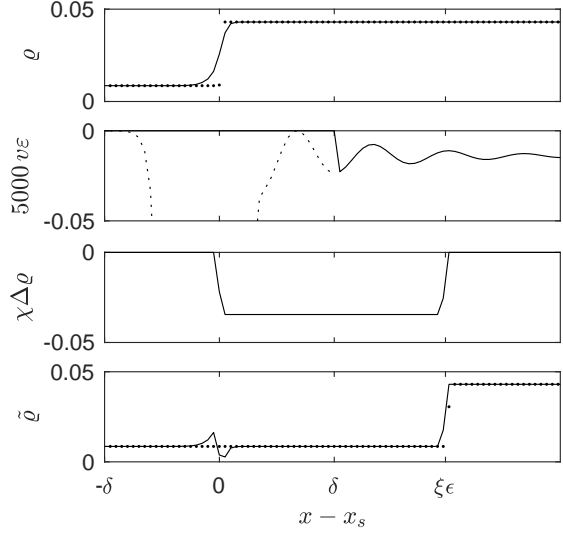}
\caption{Detailed view of the components of the calculus \eqref{eq::calculus} in the vicinity of the shock at $x_s$: Density solution of the Euler equations $\varrho$, tangential shift $v\eps$ for $|x-x_s| > \srw$, the tangential displacement of the shock $\chi\Delta \varrho$ following \eqref{eq::fv_average}, and the final approximation $\tilde{\varrho}$. For $\rho$ and $\tilde{\varrho}$, the dots denote the analytic solution in compliance with \eqref{eq::fv_average}. The dotted line complementing $v\eps$ shows the omitted data in $- \srw < x-x_s < \srw$. For this illustration the perturbation is $\eps\approx 4\mathrm{e}^{-4}$.}\label{fig::euler::calculus}
\end{figure}

The application of the calculus is illustrated in \fig{\ref{fig::euler::calculus}} for the density $\varrho$.
First, the shock at $t=t_{final}$ is shown and compared to the analytical solution illustrated by the dotted line. The shock location is accurately predicted. For $x > x_s$, the solution is not monotone, however, the under- and overshoots are six orders of magnitude smaller than the post shock density $\varrho_r$.
Compared to \fig{\ref{fig::burgers::calculus}} and considering that $\srw$ is four times greater than for the Burgers simulation, the numerical representation of the shock expands over more cells.
The approximation $\tilde{\varrho}$ does agree with the analytic solution, whereas the oscillation at the shock location of $\rho$ is more intense than for $\tilde{u}$ in \fig{\ref{fig::burgers::calculus}}. This is plausible since the steps of the discretized characteristic function $\chi$ span over one cell. The more cells the numerical representation of the shock requires, the more deviation occurs in the approximation due to the summation of $\varrho$ and $\chi\Delta\varrho$.

The tangential shift $v\eps$ in \fig{\ref{fig::euler::calculus}} is complemented by a dotted line showing the omitted data in the range $- \srw < x-x_s < \srw$.
Analogous to the derivative of $u$ in \fig{\ref{sufig::burgers::derivative_u}}, the numerical dissipation generates huge erroneous values of $v\eps$ in the vicinity of the shock at $x-x_s = 0$. Compared to the data range of the figure, the value is greater by four orders of magnitude.
At $x-x_s > \srw/2$, the tangential shift $v\eps$ undergoes oscillations that decrease with increasing $x$. In the range $\srw/2 < x < 1.5\srw$, the oscillations have the same order of magnitude than the tangential shift itself. This behavior is not present for the example of the Burgers equation and may be responsible for the divergence of the derivative of the shock location $\xi$ for values $C < 20$.
Note that the tangent $v$ enters $\xi$ by the chain rule integrating \eqref{SRH}.

\begin{figure}
\centering
\includegraphics[scale=\figscale]{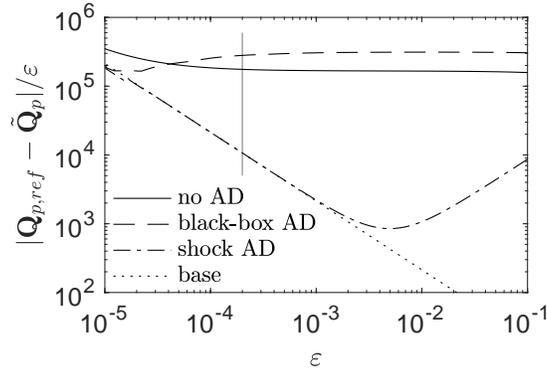}
\caption{Convergence of the tangential shift computing the propagating shock:
the tangential shift is compared to the analytic solution applying
``no AD'' no tangential shift, 
``black-box AD'' tangential shift of solution $\mathbf{Q}_p$,
``shock AD'' tangential shift using the full tangent $(v,\xi)$,
``base'' reference error at $\eps = 0 $ divided by $\eps$.
}\label{fig::euler::convergence::perturbation}
\end{figure}

Similar to \fig{\ref{sufig::burgers::convergence::perturbation}}, the convergence of the error of the approximation $\tilde{\mathbf{Q}}_p$ is analyzed in \fig{\ref{fig::euler::convergence::perturbation}}. 
Again, $\eps_{min} = 1\mathrm{e}^{-5}$ moves the shock by the cell width, $\eps_{max} = 0.1$ leads to a displacement of the shock to the end of the computational domain, and $\eps_{\dagger}$ denotes a shift by $\srw$ illustrated by the vertical line. The results show the findings analyzed in \fig{\ref{sufig::burgers::convergence::perturbation}}, i.e., the convergence of the shock AD error toward the reference error ``base'' whereas the black-box AD error does not show convergence for $\eps > \eps_{\dagger}$. Thus, the findings emphasize the validity of the theoretical calculus and numerical implementation for the Euler equations. In brief, the {shock AD} method provides correct sensitivities and the {black-box AD} application fails due to the occurrence of the discontinuity.
\section{Conclusion}\label{sec::conclusion}
Numerical computations of sensitivities of flows governed by hyperbolic 
equations remain  challenging  due to the non-existence of
a differential in any $L^p-$space. Hence, applying black--box algorithmic differentiation 
to a numerical scheme is likely to fail. The concept of tangent vectors has been 
introduced  \cite{BressanMarson1995ab} to provide an analytical framework 
for a suitable differential of such flows. In this paper, we have shown how this concept 
may be included within  finite--volume  schemes focusing in particular on the required
extension and modification necessary to apply algorithmic differentiation. In the spatially
one--dimensional case, we illustrate that the proposed algorithm leads
to suitable approximations for the sensitivity of hyperbolic flows. Results
have been presented for the Burgers equation as well as the Euler equations to 
highlight the applicability of the  algorithmic method.

\section*{Acknowledgments}
This work has been supported by DFG HE5386/18,19,  DFG 320021702/GRK2326 as well as ERS Seedfund 'Predicitive Hierarchical Simulation' 
of RWTH Aachen University.

\appendix

%
%
%
%
%
%
%

\section{Appendix on definitions and supplementary lemmas }\label{appendix}

In this section we collect definition and statements of reference \cite{BressanMarson1995aa,Bressan2000aa}. They are given for sake of completeness.

\begin{defn}[Continuous path]
	A mapping $\gamma:[a,b] \to L^1(\mathbb{R}^n)$ is called a continuous path, if 
	$\gamma$ is continuous on the interval $[a,b]$ with respect to $L^1-$norm, i.e., 
	$$\forall x \in [a,b]: \; \lim_{\eps \to 0} \| \gamma(x+\eps) - \gamma(x) \|_{L^1} =0.$$
\end{defn}

\begin{defn}[Broad solution]\label{broad solution}
	Consider the quasi--linear partial differential equation 
	\begin{equation}\label{broad sol}
	u_t(t,x) + A(t,x) u_x(t,x) = h(t,x,u),
	\end{equation}
	where $A \in \mathbb{R}^{n\times n}$ is strictly hyperbolic, Lipschitz and $h$ is measurable w.r.t. $(t,x)$ and Lipschitz 
	continuous w.r.t. $u.$ Assume  an initial condition $u(0,x) = u^0(x)$ with $u_0 \in L^1(\mathbb{R};\mathbb{R}^n)$.
	Denote by $\ell_i, r_i$  the $i$th left and right eigenvectors
	of $A.$ Denote by $\lambda_i$ the $i$th eigenvalues of $A.$ We denote by $t\to y_i(t;\tau,\xi)$ 
	the solution to the Cauchy problem 
	$$ \frac{d}{dt} y(t) = \lambda_i(t,y(t)), \; y(\tau)=\xi.$$
	Denote by $<,>$ the scalar product on $\mathbb{R}^n$ and by  
	$$g_i := < \ell_i, h> + < \partial_t \ell_{i} + \lambda_i \partial_x \ell_i, u>, \; u = \sum u_i r_i. $$
	We define a broad solution $u=\sum u_i r_i$ to equation \eqref{broad sol} 
	as a locally integrable function fulfilling 
	$$ \frac{d}{dt} u_i(t,y_i(t;\tau,\xi)) = g_i\left(t,y_i(t;\tau,\xi), u(t,y_i(t;\tau,\xi)\right)$$
	in the sense that for a.e. $(\tau,\xi)$ and all $i=1,\dots,n$ the following holds
	$$u_i(\tau,\xi) = u^0_i(y_i(0;\tau,\xi)) + \int_0^\tau 
	g_i\left(s,y_i(s;\tau,\xi), u(s,y_i(s;\tau,\xi)\right) ds.$$
\end{defn}

The main result used in this work 
is  \cite[Theorem 2.2]{BressanMarson1995aa}. We recall the statement for convenience. 

Consider the equation 
\begin{equation}\label{app1}
\partial_t u + \partial_x F(u) = h(t,x,u).
\end{equation}
supplemented with initial data $u(0,x)=u_0(x)$ and 
the assumptions
\begin{itemize}
	\item[(H1)] The vector field $F:\Omega \to \mathbb{R}^n$ is $\mathcal{C}^2$
	where $\Omega \subset \mathbb{R}^n$ is closed and bounded. For each $u \in \Omega$ the matrix $A(u)=DF(u)$ has $n$ real distinct eigenvalues. Its eigenvalues $\lambda_i$ and its left and right eigenvectors $\ell_i$ and $r_i,$ respectively, are normalized such that $<\ell_i,r_j>=\delta_{ij}$. Denote by 
	$$A(u,v) = \int_0^{1} A( \theta u + (1-\theta) v) d\theta$$
	with corresponding eigenvectors $\ell_i(u,v), \; r_i(u,v)$ and eigenvalues $\lambda_i(u,v).$ 
	Suppose that 
	$\ell_i(u,v), r_i(u,v)$ and $\lambda_i(u,v)$ are uniformly bounded for all $u,v\in \Omega.$ 
	
	\item[(H2)] Denote by $\hat{\lambda}$ the uniform bound on $\lambda_i(i,v)$ for all $i.$ 
	Then, solutions to \eqref{app1} are considered in the domain 
	$$\mathcal{D}:=\{ (t,x): 0\leq t \leq T, x \in [a+\hat{\lambda} t, b-\hat{\lambda} t] \}$$
	Assume further that the function $h:\mathcal{D}\times\Omega \to \mathbb{R}^n$ is bounded and continuously differentiable.
	
	\item[(H3)] Whenever $u^+ \in \Omega$ and $u^- \in \Omega$ are connected
	by  a shock or a contact discontinuity, say of the $k$th characteristic family, the 
	linear system 
	\begin{align*}
	0 = \Phi_i( u^+, u^-, w^+, w^-)  =\sum\limits_{j=1}^n 
	< D\ell_i(u^+,u^-) \cdot( w^+_j r_j^+, w_j^- r_j^-), u^+-u^- >  + \\
	\sum\limits_{j=1}^n < \ell_i(u^+,u^-), w_j^+ r_j^+ - w_j^- r_j^->,  \; \forall 
	i \not = k
	\end{align*}
	can be uniquely solved in terms of the outgoing variables $w^\pm_{j^\pm}$ 
	$j^{\pm} \in \{ j^- : j <k \} \cup \{ j^+: j>k \}=:\mathcal{O}.$ 
	Assume that the function $W_j$ defined by  
	$$w_{j}^\pm = W_{j^{\pm}}(u^+,u^-)( (w_j)_{ j^\pm \not \in \mathcal{O}  } ), \; j \not = k, j^\pm \in \mathcal{O} $$   
	satisfies a bound of the form 
	$$ \| W_{j^{\pm}}(u^+,u^-)( (w_j)_{ j^\pm \not \in \mathcal{O}  } ) \| \leq C 
	\| (w_j)_{ j^\pm \not \in \mathcal{O}  }  \|$$
\end{itemize}
Here, $r_j^\pm = r_j(u^\pm)$.  For a definition of the class of functions which are piecewise Lipschitz with 
simple discontinuities we refer to \cite{BressanMarson1995aa}.

\begin{thm}\label{main theorem}
	Let the assumptions $(H1)-(H3)$ hold true. Let $u$ be a piecewise Lipschitz 
	continuous solution to equation \eqref{app1} with $u^0$ in the class PLSD.
	Let $(v_0,\xi_0)\in L^1 \times \mathbb{R}^N$ be a tangent vector to $u^0$ generated
	by a regular variation $\gamma:\delta \to u^0_\delta,$ Let $u_\delta$ be the solution of equation \eqref{app1} 
	with initial condition $u^0_\delta.$ Then, there exists $\tau_0>0$
	such that for all $t \in [0,\tau_0]$ the path $\bar{\gamma}:\delta \to u_\delta$ 
	is a regular variation for $u_\delta(t,\cdot)$ generating the tangent vector
	$(v(t),\xi(t))\in L^1 \times \mathbb{R}^N.$ The vector is the unique broad solution 
	of the initial boundary value problem 
	\begin{align*}
	\xi(0) = \xi_0, v(0,x)=v_0(x), \\
	v_t + A(u) v_x + ( DA(u)v )u_x = h_u(t,x,u)v, 
	\end{align*} 
	outside the discontinuities of $u$ while for $\alpha=1,\dots,N$ 
	\begin{align*}
	< D\ell_i(u^+,u^-) \cdot (v^+ + \xi_\alpha u_x^+, v^- + \xi_\alpha u_x^-), u^+-u^-> = \\
	+ < \ell_i(u^+,u^-), v^+ + \xi_\alpha u_x^+ - v^- - \xi_\alpha u_x^- >, i \not = k_\alpha, \\
	\frac{d}{dt} \xi_\alpha = D\lambda_{k_\alpha} (u^+,u^-) ( v^+ + \xi_\alpha u_x^+, v^- + \xi_\alpha u_x^- )
	\end{align*}
	along each line $x=x_\alpha(t)$ where $u$ suffers a discontinuity in the $k_\alpha$  characteristic direction.
\end{thm}

The technical details are given in \cite{BressanMarson1995aa}. 

\end{document}